\documentclass{article}
\usepackage{xcolor}
\usepackage{mathrsfs}
\usepackage{amssymb}
\usepackage{amsmath}
\usepackage{multirow}
\usepackage{bbm}
\usepackage{csquotes}
\usepackage{soul}
\usepackage{natbib}
\usepackage{graphics}
\usepackage{graphicx}
\usepackage[a4paper, total={6in, 8in}]{geometry}

\usepackage{float}
\usepackage[hidelinks]{hyperref}
\usepackage[ruled]{algorithm2e}
\usepackage{doi}

\numberwithin{equation}{section}
\DeclareMathOperator*{\sign}{sign}
\def\ud{\, \mathrm{d}}

\newtheorem{theorem}{Theorem}[section]

\newtheorem{Definition}[theorem]{Definition}

\begin{document}
\title{Investigation and Development of the Methodologies for Simulating Self-similar Processes}

\author{Qidi Peng\footnote{Corresponding author. Institute of Mathematical Sciences, Claremont Graduate University. Email: qidi.peng@cgu.edu.} and  William Wu\footnote{Rancho Cucamonga High School. Email: wwu49387@gmail.com}}
\date{}
\maketitle

\begin{abstract}
This paper is devoted to the study of simulating a large class of self-similar processes. Since most current simulation approaches are limited to case-by-case studies,  every existing approach has its constraints and flaws; hence a general and efficient simulation approach is in demand. Our study sheds some light in this direction. The paper's contributions are bi-fold. First, reviews and improvements are made to some existing methods for simulating specific self-similar processes. Second, we propose a novel method to simulate a general self-similar process, where we use a modified inverse Lamperti transformation to transform self-similarity to stationarity. Successful applications are made to simulate fractional Brownian motion and sub-fractional Brownian motion.

\begin{flushleft}

\end{flushleft}
\end{abstract}

\section{Introduction}
\label{sec:introduction:selfsimilar}
Self-similar processes are an essential and large class of stochastic processes that exhibit the phenomenon of \enquote{self-similarity}. That means part of the process's path behaves similarly to its entire path. Self-similarity is featured in many classical stochastic processes in probability theory, such as Brownian motion, fractional Brownian motion \cite{Mandelbrot1968}, linear fractional stable motion \cite{Taqqu1994}, fractional Poisson process \cite{laskin2003fractional}, and Rosenblatt process \cite{dobrushin1979non,taqqu1979convergence}. Some random fields and random sheets also exhibit self-similarity, with examples such as fractional Brownian fields \cite{dobric2006fractional}, linear fractional stable fields and sheets \cite{ayache2007local,ding2023linear}, etc. It should be noted that self-similar processes can be naturally extended to locally asymptotically self-similar processes \cite{Boufoussi2008}, with  multifractional Brownian motion and linear multifractional stable motion being pragmatic examples   \cite{vehel1998introduction,Benassi1997Elliptic,stoev2004simulation,ding2023linear}. Modern applications of self-similar processes include modeling instrument values in finance, volume of flows in hydrology, weather conditions in meteorology, atomic movement in material science and physics, etc.  \cite{rostek2013note,haltas2011scale,elsaesser2010observed,sun2018self}. 
In the literature, \enquote{simulation of self-similar process}  has received wide attention. For example, simulations of fractional ARIMA \cite{stoev2004simulation}, fractional Brownian motion \cite{embrechts2009selfsimilar},  linear fractional stable motion \cite{stoev2004simulation}, and self-similar teletraffic \cite{jeong2002modelling} have been heavily studied. So far, the methods for simulating self-similar processes work case by case and each has its flaws and constraints. Indeed, \cite{embrechts2009selfsimilar} in Chapter 7 stated that  ``very little specific tools are available`` for simulating a general self-similar process. 

 In the following sections, we list some main existing simulation methodologies for self-similar processes.  We explain how the above methods are applied to simulate particular self-similar processes and discuss their pros and cons. Furthermore, we develop a new method to overcome the issues displayed by the above approaches. Therefore, the goals of this paper are the following:
 \begin{description}
 \item[(i)] We investigate and summarize simulation methodologies existing in the literature for a number of well-known self-similar processes. We provide pseudocodes and sources of implementations of these simulation algorithms, and discuss their pros and cons.\\
 \item[(ii)] We develop a novel method for simulating self-similar Gaussian  process based on a modified inverse Lamperti transformation. This method has the potential to be extended for simulating more general self-similar processes.
 \end{description}
 Through this paper, we classify self-similar processes into stationary increments Gaussian processes, processes with integral representations, and processes with wavelet decompositions.

 The paper is organized as follows. In Sections \ref{Sec2} - \ref{Sec5}, we investigate, respectively, the existing simulation approaches for stationary increments self-similar processes, multifractional Brownian motion, self-similar processes with non-stationary increments, and linear fractional stable motion. These examples cover self-similar processes from Gaussian processes to non-Gaussian ones and from stationary increments processes to  non-stationary increments ones. Moreover, we discuss simulation of locally asymptotically  self-similar processes, which are natural extensions of self-similar processes. Our significant contribution is in Section \ref{Sec6}, where a brand new simulation algorithm is introduced to simulate Gaussian self-similar processes. This method is based on a modified inverse Lamperti transformation and it has the potential to be applied to simulate a large family of self-similar processes. We take the simulation of fractional Brownian motion and sub-fractional Brownian motion as two examples. In Section \ref{Sec7}, we conclude.

 \section{Simulation of Stationary Increments Self-similar Gaussian Processes}
 \label{Sec2}
 The self-similar process considered in this paper is defined below.
\begin{Definition}
\label{def:selfsimilar}
We say $X=\{X(t)\}_{t\in \mathbb{R}}$ is a self-similar process with self-similarity index $H\in(0,1)$, if
\begin{equation}
\label{self-similar}
    \left\{X(at)\right\}_{t\in\mathbb R}\stackrel{\mbox{f.d.d.}}{=}\left\{|a|^HX(t)\right\}_{t\in\mathbb R}, ~\mbox{for any}~ a\ne0,
\end{equation}
where $\stackrel{\mbox{f.d.d.}}{=}$ denotes equality in finite dimensional distribution: for two stochastic processes $\{X(t)\}_{t\in\mathbb R}$ and $\{Y(t)\}_{t\in\mathbb R}$, $\{X(t)\}_{t\in\mathbb R}\stackrel{\mbox{f.d.d.}}{=}\{Y(t)\}_{t\in\mathbb R}$ means 
$
(X(t_i))_{i=1,
\ldots,n}\stackrel{\mbox{law}}{=}(Y(t_i))_{i=1,
\ldots,n}
$
for any $n\ge1$ and any $t_1,\ldots,t_n\in\mathbb R$.
\end{Definition}
 Let $\{X(t)\}_{t\ge0}$ be a self-similar process. By self-similarity, it suffices to simulate the sample path of $\{X(t)\}_{t\in[0,1]}$ in order to obtain any sample paths of $\{X(t)\}_{t\in[a,b]}$ for $0<a<b$, because (see (\ref{self-similar}))
 $$
 \left\{X(t)\right\}_{t\in[a,b]}\stackrel{\mbox{f.d.d.}}{=}\left\{b^HX(t)\right\}_{t\in[a/b,1]}.
 $$
 Therefore, in the rest of the paper, we only consider simulating a discretized path of $\{X(t)\}_{t\in[0,1]}$: $\left(X(1/N),~X(2/N),~\ldots,X((N-1)/N),~X(1)\right)$, where $N$ denotes the number of nodes in the trajectory.  The major techniques in play are Cholesky's decomposition and fast Fourier transformation. The simulation approaches introduced in this section apply to Brownian motion, fractional Brownian motion, fractional Gaussian noise, fractional ARIMA, and fractional teletraffic data process.  

 Recall that a fractional Brownian motion $\{B^H(t)\}_{t\in\mathbb R}$ with Hurst parameter $H\in(0,1)$ \cite{kolmogorov1940wienersche,Mandelbrot1968} is a zero-mean Gaussian stochastic process with almost surely $B^H(0)=0$, continuous trajectories, and its covariance function  is given by: for $s,t\in\mathbb R$,
\begin{equation}
\label{cov_fBm}
\mathbb Cov\left(B^H(s),B^H(t)\right)=\frac{1}{2}\left(|t|^{2H}+|s|^{2H}-|t-s|^{2H}\right).
\end{equation}
Note that when $H=1/2$, the fractional Brownian motion becomes a standard Brownian motion, denoted by $\{B(t)\}_{t\in\mathbb R}$. Since  $\{B(t)\}_{t\in[0,1]}$ is a L\'evy process (it has independent stationary increments), it can be decomposed into the cumulative sum of Gaussian white noise:
\begin{equation}
\label{sim:BM}
\left\{\begin{array}{lll}
&B(0)=0~\mbox{a.s.};&\\
&B\left(\frac{k}{N}\right)=\sum_{j=0}^{k- 1}\varepsilon_{j,N},&~\mbox{for}~k=1,\ldots,N,
\end{array}\right.
\end{equation}
where $\{\varepsilon_{j,N}\}_{j\in\{0,\ldots,N-1\}}=\left\{B\left((j+1)/N\right)-B\left(j/N\right)\right\}_{j\in\{0,\ldots,N-1\}}$ is a sequence of i.i.d.  zero-mean Gaussian variables with $\mathbb Var(\varepsilon_{j,N})=1/N$. (\ref{sim:BM}) yields the simulation of $(B(1/N),\ldots,B(1))$. 

For simulating discretized trajectories of fractional Brownian motion $\{B^H(t)\}_{t\in[0,1]}$ with $H\ne 1/2$, we  face the challenge of simulating correlated Gaussian random vectors. Multiple approaches have been suggested to address this challenge.
\begin{description} 
\item[Cholesky's decomposition-based method : ] This method is often used in a general approach to simulate a Gaussian vector (process). Let $\mathbb N=\{0,1,2,\ldots\}.$ Assume that a discretized Gaussian process  $\{G(k)\}_{k\in\mathbb N}$ has covariance function
$$
\gamma_{i,j}=\mathbb Cov(G(i),G(j)),~\mbox{for}~i,j\in\mathbb N.
$$
Then, for any time index $n\in\mathbb N$, given $(G(0),\ldots,G(n))$, the conditional distribution of $G(n+1)$ is also Gaussian with its mean and variance specified below: 
\begin{eqnarray}
\label{sim_Gaussian}
&&G(n+1)|(G(0),\ldots,G(n))\nonumber\\
&&\sim \mathcal N\left(
C(n)^TM(n)^{-1}
\begin{pmatrix}
G(0)\\
\vdots\\
G(n)
\end{pmatrix},~\gamma_{n+1,n+1}-C(n)^TM(n)^{-1}
C(n)\right),
\end{eqnarray}
where 
\begin{itemize}
\item The column matrix
$
C(n)=\begin{pmatrix}
\gamma_{n+1,0}\\
\vdots\\
\gamma_{n+1,n}
\end{pmatrix}.
$
\item $M(n)$ is the covariance matrix of $(G(0),\ldots,G(n))$, it  satisfies the following recursion: for $n\ge0$,
$$
M(n+1)=
\begin{pmatrix}
M(n)&C(n)\\
C(n)^T&\gamma_{n+1,n+1}
\end{pmatrix}.
$$
\end{itemize}
The general method for simulating a Gaussian process is then: 
\begin{equation}
\label{sim:G}
\left\{\begin{array}{lll}
&\mbox{Set}~G(0)=x_0,~\mbox{where $x_0\in\mathbb R$ is a given initial value};&\\
&\mbox{For}~k=1,\ldots,n,~G\left(k\right)~\mbox{is simulated via the distribution in (\ref{sim_Gaussian}).}
\end{array}\right.
\end{equation}
Note that one has to compute $C(n)^TM(n)^{-1}$ for each $n\ge1$ in this simulation method. Several Cholesky's decomposition-based implementation methods are therefore suggested to compute this term, see \cite{michna1998self,michna1999tail,michna2000ruin,embrechts2009selfsimilar}. We point out that this general simulation approach applies to the simulation of fractional Brownian motion with covariance function (\ref{cov_fBm}).\\
\item[Fast Fourier transformation-based method : ] Since fractional Brownian motion is a stationary increments Gaussian process, we can alternatively apply Davies and Harte's method \cite{davies1987tests} for simulating it. The implementation relies on the fast Fourier transformation. Let the auto-covariance function of the increments of fractional Brownian motion  be given by
$$
\gamma_{H}\left(\frac{k}{n}\right)=\mathbb Cov\left(B^H\left(\frac{k+1}{n}\right)-B^H\left(\frac{k}{n}\right),~B^H\left(\frac{1}{n}\right)\right).
$$
We adjust Davies and Harte's \cite{davies1987tests} algorithm to simulate the fractional Brownian motion $\{B^H(t)\}_{t\in[0,1]}$ below:\\
\begin{minipage}{16cm}
\begin{algorithm}[H]
\label{algo:sim_fBm_FFT}
\caption{Simulation of fractional Brownian motion via the fast Fourier transform.}
\LinesNumbered 
\KwIn{\textbf{\textbf{\textit{path length}} $n$; \textit{auto-covariance function}} $\gamma_H$.}
\For{$k= 1,\ldots,2n-2$}{
$\alpha_k \longleftarrow \frac{2\pi(k-1)}{2n-2}$\;
$f_k \longleftarrow \max\left\{\mbox{Re}\left(\sum\limits_{j=1}^{n-1}\gamma_{H}\left(\frac{j-1}{n}\right)e^{i(j-1)\alpha_k/n}+\sum\limits_{j=n}^{2n-2}\gamma_H\left(\frac{2n-j-1}{n}\right)e^{i(j-1)\alpha_k/n}\right),0\right\}$\;}
$V_1 \longleftarrow0$\;
$V_n \longleftarrow0$\;
\textbf{\textit{Simulate independent zero-mean Gaussian variables $(U_1,\ldots,U_n)$ and $(V_2,\ldots,V_{n-1})$  with}}\\
$\mathbb Var(U_1)=\mathbb Var(U_n)=2$\;
$\mathbb Var(U_k)=\mathbb Var(V_k)=1$~for $k\ne 1,n$\;
\For{$k = 1,\ldots,n$}{
$Z_k \longleftarrow U_k+iV_k$\;}
\For{$k = n+1,\ldots,2n-2$}{
$Z_k \longleftarrow U_{2n-k}-iV_{2n-k}$\;}
\For{$t = 1,\ldots,n$}{
$Y_t \longleftarrow \frac{1}{2\sqrt{n-1}}\sum\limits_{k=1}^{2n-2}\sqrt{f_k}e^{i(t-1)\alpha_k/n}Z_k$\;}
\For{$t = 1,\ldots,n$}{
$B^H\left(\frac{t}{n}\right) \longleftarrow \sum\limits_{k=1}^tY_k$\;}
\KwOut{ $(B^H(1/n),~B^H(2/n),\ldots,B^H(1))$.}
\end{algorithm}
\end{minipage}
The notation $\mbox{Re}(\bullet)$ denotes the real part of a complex value. The advantage of this method is that one can use the fast Fourier transform to efficiently compute $f_k$ and $Y_t$. Moreover, Line 3 guarantees that $f_k\ge0$, because the square root operator of $f_k$ will be needed in the computation of $Y_t$ in Line 17. Because of inherent computational error, it happens that $f_k=0$ for some $k$, which subsequently underestimates the covariances of the increments of fractional Brownian motion. Therefore, this fast Fourier transform method carries the risk of destroying the long-memory property of fractional Brownian motion. This issue has been mentioned in \cite{asmussen1998stochastic,embrechts2009selfsimilar} and is caused by the fact that the positivity condition for $f_k$ is not always satisfied. Inspired by the simulation procedure in  \cite{davies1987tests}, \cite{WoodChan1994} have proposed the circulant embedding approach for simulating stationary Gaussian random fields. In this latter framework, the positivity condition issue is addressed  through \enquote{treating the dimension of the circulant as an integer parameter that is chosen to ensure that the circulant is positive definite}. To elaborate, Algorithm \ref{algo:sim_fBm_FFT} applies the fast Fourier transform on the so-called Toeplitz matrix
$$
T=\begin{pmatrix}
\gamma_H(0)&\gamma_H\left(\frac{1}{n}\right)&\cdots&\gamma_H\left(\frac{n-1}{n}\right)\\
&&&\\
\gamma_H\left(\frac{1}{n}\right)&\gamma_H\left(0\right)&\cdots&\gamma_H\left(\frac{n-2}{n}\right)\\
&&&\\
\vdots&\vdots&\ddots&\vdots\\
&&&\\
\gamma_H\left(\frac{n-1}{n}\right)&\gamma_H\left(\frac{n-2}{n}\right)&\cdots&\gamma_H\left(0\right)
\end{pmatrix}.
$$
Note that $T$ is not necessarily non-negative definite. Therefore,  \cite{WoodChan1994} suggest to embed $T$ in a circulant covariance matrix
$$
C=\begin{pmatrix}
c_0&c_1&\cdots&c_{m-1}\\
&&&\\
c_{m-1}&c_0&\cdots&c_{m-2}\\
&&&\\
\vdots&\vdots&\ddots&\vdots\\
&&&\\
c_1&c_2&\cdots&c_0
\end{pmatrix},
$$
where $m\ge 2(n-1)$ is some integer and 
$$
c_j=\left\{
\begin{array}{ll}
\gamma_H\left(\frac{j}{m}\right)&~\mbox{if}~0\le j\le \frac{m}{2};\\
\gamma_H\left(\frac{m-j}{m}\right)&\mbox{if}~\frac{m}{2}<j\le m-1.
\end{array}\right.
$$
Next, \cite{WoodChan1994} suggest either to increase $m$ or use the non-negative definite part of $C$ in order to match the positivity condition.\\
\item[Moving average stochastic integral representation method : ] Recall that the fractional Brownian motion can be equivalently defined via the moving average integral representation \cite{taqqu1979self}:
\begin{equation}
\label{def:fBm_ma}
B^H(t)=C_H\int_{-\infty}^t\left((t-u)^{H-1/2}-\max(-u,0)^{H-1/2}\right)\ud B(u),
\end{equation}
where the constant $C_H$ is chosen such that $\mathbb Var(B^H(1))=1$:
$$
C_H=\left\{\int_{-\infty}^0\left((t-u)^{H-1/2}-(-u)^{H-1/2}\right)^2\ud u+\frac{1}{2H}\right\}^{-1/2}.
$$
Based on (\ref{def:fBm_ma}), one can use the truncation method to simulate a discrete path of fractional Brownian motion. However, since the stochastic integral (\ref{def:fBm_ma}) is indefinite, the truncation might also destroy the long-memory property \cite{embrechts2009selfsimilar}. \\
\item[Series representation methods : ] There exist other simulation methods based on representation results of fractional Brownian motion. For example,  the representation of fractional Brownian motion in Theorem 5.2 given in \cite{norros1999elementary} does not require  truncation. Therefore, the long-memory property is preserved by the simulation method based on the following integral representation:
$$
B^H(t)=\int_0^tz(t,s)\ud W(s),
$$
where $\{W(s)\}_{s\ge0}$ is a standard Brownian motion and 
\begin{eqnarray*}
&&z(t,z)=\left(\frac{2H\Gamma(3/2-H)}{\Gamma(H+1/2)\Gamma(2-2H)}\right)^{1/2}\nonumber\\
&&\times\left(\left(\frac{t}{s}\right)^{H-1/2}(t-s)^{H-1/2}-(H-\frac{1}{2})s^{1/2-H}\int_s^tu^{H-3/2}(u-s)^{H-1/2}\ud u\right).
\end{eqnarray*}
For fractional Brownian motion with the long-memory property ($H>1/2$), we can simplify the above expression to
$$
z(t,s)=\left(\frac{2H\Gamma(3/2-H)}{\Gamma(H+1/2)\Gamma(2-2H)}\right)^{1/2}{}_2F_1\left(\frac{1}{2}-H,H-\frac{1}{2},H+\frac{1}{2},1-\frac{t}{s}\right),
$$
where ${}_2F_1$ is the Gauss hypergeometric function. Finally, there exist alternative simulation methods based on wavelet representations. For example, a mean-square sense wavelet representation of fractional Brownian motion is given in \cite{flandrin1992wavelet}, where, for any given resolution $2^J$ with $J\in\mathbb Z$,
$$
B^H(t)=2^{-J/2}\sum_{n=-\infty}^{+\infty}a_{J,n}\phi(2^{-J}t-n)+\sum_{j=-\infty}^J2^{-j/2}\sum_{n=-\infty}^{+\infty}d_{j,n}\psi(2^{-j}t-n),
$$
where $\psi(t)$ is the basic wavelet function and $\phi(t)$ is the scaling function associated to $\psi(t)$ \cite{mallat1989theory}; the approximation coefficients are
\begin{eqnarray*}
&&a_{j,n}=2^{-j/2}\int_{-\infty}^{+\infty}B^H(t)\phi(2^{-j}t-n)\ud t;\\
&&d_{j,n}=2^{-j/2}\int_{-\infty}^{+\infty}B^H(t)\psi(2^{-j}t-n)\ud t.
\end{eqnarray*}
There exist other wavelet representations of fractional Brownian motion, for example, in \cite{whitcher2001simulating,albeverio2012fractional,coeurjolly2013wavelet}. The wavelet decomposition methods may suffer from accuracy, high complexity issues, and the destruction of the long-memory property.
\end{description}
In conclusion, so far, the most promising method for simulating fractional Brownian motion stems from \cite{WoodChan1994}, thanks to the fact that the method does not destroy the long-memory property and is relatively accurate and efficient.

For simulating other self-similar Gaussian processes very similar to fractional Brownian motion, such as fractional Gaussian noise and fractional ARIMA, we refer the readers to \cite{beran2017statistics,doukhan2002theory}. Murad Taqqu's website \url{http://math.bu.edu/people/murad} offers great insight on simulation strategies for self-similar processes \cite{embrechts2009selfsimilar}.

Some real-world data also display self-similar and light-tail properties, such as fractional teletraffic data. Simulating these types of data plays an important role in both industrial settings and academic research. Various methods have been suggested to simulate them and are more or less inspired by those  for simulating fractional Brownian motion; see the FGN-DW methods in \cite{jeong1999fast,jeong2002modelling,jeong2003generation}. 

\section{Simulation of Multifractional Brownian Motion}
Multifractional Brownian motion, as one of the most natural extensions of fractional Brownian motion, is a continuous-time Gaussian process with a time-dependent self-similarity parameter. It is not self-similar but locally asymptotically self-similar \cite{Benassi1997Elliptic,Boufoussi2008}. Unlike fractional Brownian motion, multifractional Brownian motion has non-stationary increments when its Hurst functional parameter is not a constant. In the literature, multifractional Brownian motion has been introduced separately and independently by both  \cite{vehel1998introduction} and  \cite{Benassi1997Elliptic}. In \cite{vehel1998introduction}, multifractional Brownian motion can be defined via the moving-average representation
\begin{equation}
\label{def:mBm_ma}
X^{(1)}(t)=\int_{-\infty}^{+\infty}\left((t-u)_+^{H(t)-1/2}-(-u)_+^{H(t)-1/2}\right)\ud W(u),
\end{equation}
where we denote $(x)_+=\max\{x,0\}$ for any real number $x$, the function $H$ is continuous, taking values in $(0,1)$, and $\ud W(u)$ denotes an independently scattered standard Gaussian measure on $\mathbb R$. In \cite{Benassi1997Elliptic}, multifractional Brownian motion is defined via the so-called harmonizable integral representation as
\begin{equation}
\label{def:mBm_ha}
X^{(2)}(t)=\int_{-\infty}^{+\infty}\frac{e^{it\xi}-1}{|\xi|^{H(t)+1/2}}\ud\widetilde W(\xi),
\end{equation}
where $\ud\widetilde W(\xi)$ denotes a complex-valued Gaussian measure. Note that the two definitions (\ref{def:mBm_ma}) and (\ref{def:mBm_ha}) have different covariance structures  due to different scaling coefficients \cite{Stoev2006}. In this paper, we discuss only simulation methods of the multifractional Brownian motion $X(t)$ defined by \cite{chan1998simulation}, which is a modification of $X^{(1)}(t)$. Simulations of both $X^{(1)}(t)$ and
 $X^{(2)}(t)$ can be obtained from simple adjustments. 

Recall that the version of the multifractional Brownian motion considered in \cite{chan1998simulation} is the following:
 $$
 X(t)=\widetilde C(H)\int_{-\infty}^{+\infty}\left(1-\cos(st)-\sin(st)\right)|s|^{-H-1/2}\ud W(s),
 $$
 where $\widetilde C(H)$ is a real-valued function such that $\mathbb Var(X(t))=|t|^{2H}$, and $W$ is a standard Brownian motion. $X(t)$ has the covariance function: for $s,t\ge0$,
 \begin{equation}
 \label{cov:mBm_Wood_Chan}
 \mathbb Cov(X(s),X(t))=\frac{g(H(s),H(t))}{2}\left(|s|^{H(s)+H(t)}+|t|^{H(s)+H(t)}-|t-s|^{H(s)+H(t)}\right),
 \end{equation}
 where for $H_1,H_2,H\in(0,1)$, $g(H_1,H_2)=(I(H_1)I(H_2))^{-1/2}I((H_1+H_2)/2)$ and
 $$
 I(H)=\left\{
 \begin{array}{ll}
 \frac{\Gamma(1-2H)}{H}\sin\left(\frac{(1-2H)\pi}{2}\right)&~\mbox{for}~H\in(0,1/2);\\
 \pi&~\mbox{for}~H=1/2;\\
 \frac{\Gamma(2(1-H))}{H(2H-1)}\sin\left(\frac{(2H-1)\pi}{2}\right)&~\mbox{for}~H\in(1/2,1).
 \end{array}\right.
 $$

The main simulation approaches provided by the literature consist of Wood-Chan's method \cite{chan1998simulation} and wavelet-based simulations \cite{ayache2011multiparameter,Jin2018}.

We summarize Wood-Chan's algorithm \cite{chan1998simulation} to simulate the multifractional Brownian motion $\{X(t)\}_{t\in[0,1]}$ below:\\
\begin{minipage}{16cm}
\begin{algorithm}[H]\footnotesize
\label{algo:sim_mBm_FFT}
\caption{Simulation of multifractional Brownian motion via FFT.}
\LinesNumbered 
\KwIn{\textbf{\textbf{\textit{the path length}} $n$; \textit{the discretized Hurst functional parameter}} $H_{1/(m+1)},\ldots,H_{m/(m+1)}$.}
\For{$u= 1,\ldots,m$}{
\For{$j=1,\ldots,n$}{
\textbf{\textit{Simulate $Y_{j,u}$ using the circulant embedding approach for stationary vector-valued Gaussian processes via the covariance function (see (\ref{cov:mBm_Wood_Chan})):
$$
\mathbb Cov\left(Y_{j,u},Y_{k,v}\right)=\frac{g(H_u,H_v)}{2}\left(|j-k-1|^{H_u+H_v}+|j-k+1|^{H_u+H_v}-2|j-k|^{H_u+H_v}\right);
$$}}
}}
\For{$u= 1,\ldots,m$}{
\For{$j=1,\ldots,n$}{
$Z_{j,u} \longleftarrow n^{H_u}\sum\limits_{k=1}^jY_{k,u}$\;
}}
\For{$j=1,\ldots,n$}{
\textbf{\textit{Simulate $X_{j/n}$ using some form of kriging based on the observations $\{Y_{j,u}\}_{j\in\{1,\ldots,n\},u\in\{1,\ldots,m\}}$:
\begin{eqnarray*}
&&N_j\longleftarrow\{j-q,\ldots,j+q\}\times\{[mH_{j/n}]-q,\ldots,[mH_{j/n}]+q\}\cap\{1,\ldots,n\}\times\{1,\ldots,m\};\\
&&X(j/n)\longleftarrow\sum_{(k,v)\in N_j}\mathbb Cov((Z_{k,v})_{(k,v)\in N_j})^{-1}\mathbb Cov\left((Z_{k,v})_{(k,v)\in N_j},Z_{j,j}\right)Z_{k,v};
\end{eqnarray*}}}
}
\KwOut{A sample path $(X(1/n),~X(2/n),\ldots,X(1))$.}
\end{algorithm}
\end{minipage}
Note that in Lines 12 - 13 of Algorithm \ref{algo:sim_mBm_FFT}, each $X(j/n)$ is evaluated by minimizing the variance of $X(j/n)-\sum_{(k,v)\in N_j}\beta_{k,v}Z_{k,v}$ with respect to $\beta_{k,v}$. The optimizer is given as
$$
\widehat \beta_{k,v}=\mathbb Cov((Z_{k,v})_{(k,v)\in N_j})^{-1}\mathbb Cov\left((Z_{k,v})_{(k,v)\in N_j},Z_{j,j}\right),
$$
where the covariance matrices
$$
\mathbb Cov((Z_{k,v})_{(k,v)\in N_j})=\mathbb E\left\{(Z_{k,v})_{(k,v)\in N_j}^T(Z_{k,v})_{(k,v)\in N_j}\right\}
$$
and 
$$
\mathbb Cov\left((Z_{k,v})_{(k,v)\in N_j},Z_{j,j}\right)=\mathbb E\left\{Z_{j,j}(Z_{k,v})_{(k,v)\in N_j}\right\}.
$$
The right-hand sides of the above two equations are computed using the covariance function (\ref{cov:mBm_Wood_Chan}).

Since Algorithm \ref{algo:sim_mBm_FFT} uses Wood and Chan's circulant embedding method \cite{WoodChan1994}, it therefore avoids ruining the long-memory property of the locally asymptotical self-similarity of multifractional Brownian motion. Algorithm \ref{algo:sim_mBm_FFT} itself works, and its implementation in MATLAB can be found in FracLab (\url{https://project.inria.fr/fraclab/}). However, we have found that this implementation of the fractional Brownian random field $Y_{j,u}$, $(j,u)\in\{1,\ldots,n\}\times\{1,\ldots,m\}$ 
given in Line 3 is only approximative and does not follow the approach for simulating a Gaussian random field given in \cite{WoodChan1994}. In fact, the implementation in FracLab assumes that the vectors $(Y_{1,u},\ldots,Y_{n,u})$ are independent with respect to $u\in\{1,\ldots,m\}$, although they are not. This assumption, together with Line 12 introduces more bias than Wood-Chan's method for simulating fractional Brownian motions. So far, we have not yet found a resource that exactly implements Algorithm \ref{algo:sim_mBm_FFT}. To conclude, there is still space to improve the accuracy of Wood-Chan's method for simulating multifractional Brownian motion.

Other simulation methods include the moving-average stochastic integral representation method \cite{vehel1998introduction} and Wavelet-based decomposition methods  \cite{ayache2007wavelet,ding2023linear}. Similarly to the approaches proposed for simulating fractional Brownian motion, these approaches have pros and cons. While the moving-average stochastic integral representation method truncates the integral and destroys long-range dependency, wavelet-based methods are largely inefficient and inaccurate.

\section{Simulation of Self-similar Processes with Non-stationary Increments}

The well-known fractional Brownian motion has stationary increments. \cite{WoodChan1994} have taken advantage of this feature to simulate fractional Brownian motion and multifractional Brownian motion. However, not all self-similar processes exhibit this property. Examples include sub-fractional Brownian motion \cite{bojdecki2004sub}, bi-fractional Brownian motion \cite{houdre2003example}, and tri-fractional Brownian motion \cite{ma2013schoenberg}. In these cases, Wood and Chan's methods no longer work. We take the sub-fractional Brownian motion as an example below.

The sub-fractional Brownian motion (sfBm for short) $\{S^H(t)\}_{t\ge0}$ is a
centered Gaussian self-similar process with $S^H(0)=0$ and covariance function: for $s,t\ge0$,
\begin{equation}
\label{def:sfBm}
\mathbb Cov(S^H(t),S^H(s))=t^{2H}+s^{2H}-\frac{1}{2}((t+s)^{2H}+|t-s|^{2H}),
\end{equation}
with the parameter $H\in(0,1)$. For sub-fractional Brownian motion, the existing simulation methods consist of some variants of the moving-average integral representations \cite{morozewicz2015simulation,kuang2017asymptotic}. Like other integral representation methods, the simulation approaches in \cite{morozewicz2015simulation} and \cite{kuang2017asymptotic} suffer from truncation or approximation problems.

For bi-fractional Brownian motion and tri-fractional Brownian motion which have neither non-stationary increments nor simple integral representations, we have not found any study on the simulation of these two processes in the literature, showcasing the complexity of researching and simulating self-similar processes. In Section~\ref{Sec6}, we will shed some light on this problem by introducing a brand new simulation methodology. 

\section{Simulation of Linear Fractional Stable Motion}
\label{Sec5}
In this section, we investigate the simulation of self-similar non-Gaussian processes. A typical example is linear fractional stable motion (LFSM), introduced in \cite{Taqqu1994}. Recall that the linear fractional stable motion $\{X(t)\}_{t\in\mathbb R}$ with the tail-heaviness control parameter (also called stability parameter) $\alpha\in(0,2)$ and the Hurst parameter (also called self-similarity parameter) $H\in(0,1)$  \cite{ayache2007local} is defined as follows: on a filtered probability space $(\Omega,\mathcal F,(\mathcal F_t)_{t\ge0},\mathbb P)$,  
 \begin{equation} 
\label{LFSM_definition}
X(t) = \int_{\mathbb R} g_{H}(t,s)\ud M_{\alpha}(s),~\mbox{for}~t\in\mathbb R,
\end{equation}
where:
\begin{description}
\item[(1)] In the kernel 
\begin{equation}
\label{g}
g_{H}(t,s) = \kappa\left\{(t-s)_{+}^{H-1/\alpha} - (-s)_{+}^{H-1/\alpha}\right\},
\end{equation}
$\kappa>0$ is a normalizing constant such that the scale parameter of $X(1)$, denoted by 
 $
 \|X(1)\|_{\alpha}=(\int_{\mathbb R}|g_H(1,s)|^\alpha\ud s)^{1/\alpha},
 $ 
 equals 1; we employ the convention $0^{H-1/\alpha} = 0$ for all $H\in(0,1)$ and $\alpha\in(0,2)$.\\
\item[(2)] $M_{\alpha}$ is a strictly  $\alpha$-stable random measure on $\mathbb R^N$ with Lebesgue measure as its control measure and with $\beta(s)$ as its skewness intensity (see \cite{Taqqu1994}). That is, for every Lebesgue measurable set $E \subset \mathbb R^N$ with Lebesgue measure $\lambda(E) < +\infty$, $M_{\alpha}(E)$ is a strictly $\alpha$-stable random variable with scale parameter $\sigma=\lambda(E)^{1/\alpha}$ and skewness parameter $\beta=(1/\lambda(E))\int_{E}\beta(s)\ud s$. More precisely, the characteristic function of $M_\alpha(E)$ has the following form (see Definition 1.1.6, \cite{Taqqu1994}): for $x\in\mathbb R$,
$$
\mathbb E\left[e^{ixM_\alpha(E)}\right]=\left\{
\begin{array}{ll}
\exp\left\{-\sigma^\alpha|x|^\alpha(1-i\beta(\sign(x))\tan(\pi\alpha/2))\right\}&~\mbox{if}~\alpha\ne1;\\
\exp\left\{-\sigma|x|(1+2i\beta\sign(x)\log(|x|)/\pi)\right\}&~\mbox{if}~\alpha=1,
\end{array}\right.
$$
where 
$$
\sign(x)=\left\{
\begin{array}{ll}
1&~\mbox{if}~x>0;\\
0&~\mbox{if}~x=0;\\
-1&~\mbox{if}~x<0.
\end{array}\right.
$$
\end{description}
We also denote $M_\alpha(E)\sim S_\alpha(\sigma,\beta)$. 
If $\beta(\bullet)\equiv0$, $M_\alpha$ becomes a symmetric $\alpha$-stable random measure. In this case, $M_\alpha(E)\sim S_\alpha(\sigma,0)$, i.e., 
$$
\mathbb E\left[e^{ixM_\alpha(E)}\right]=\exp\left\{-\sigma^\alpha|x|^\alpha\right\}~\mbox{for}~\alpha\in(0,2).
$$
 The special convention $0^{H-1/\alpha}=0$ in (\ref{g}) yields that, when $H=1/\alpha$, LFSM becomes the ordinary stable sheet. When $\alpha=2$, LFSM becomes fractional Brownian motion. 

 Since LFSM is non-Gaussian, the simulation approaches based on the moments, such as Cholesky's decomposition and the circulant embedding methods, do not apply anymore. The existing simulation methods are through integral representations and the Fourier series decomposition method. The integral representation methods were introduced in \cite{stoev2004simulation}. This method directly truncates the integral in (\ref{LFSM_definition}).  \cite{wu2004simulating} approximate the integral (\ref{LFSM_definition}) based on a linear process. \cite{bierme2008fourier} use the Fourier series representation of LFSM to simulate it.

 \section{A Modified Inverse Lamperti Transformation Approach}
\label{Sec6}
In this section, we develop a novel method to simulate a zero-mean Gaussian self-similar process. Our simulation methodology involves using a modified inverse Lamperti transformation, which transforms self-similar processes to stationary ones. We will take the simulation of fractional Brownian motions and sub-fractional Brownian motions as examples, but point out that this new idea has the potential to be generalized to simulate a much larger class of self-similar processes. 

Based on the above investigation of the existing simulation approaches, we believe that working on the transformation of self-similar processes to stationary processes may point in the direction of an efficient generating algorithm of self-similar processes. Therefore, the Lamperti transformation \cite{lamperti1962semi} can be applied.  
 The Lamperti transformation of a stationary process $\{Y(t)\}_{t\ge0}$ is defined to be
$$
X(t)=t^HY(\log(t)),~\mbox{for}~t>0~\mbox{and some}~H\in(0,1).
$$
The Lamperti transformation with parameter $H$ transforms a stationary process to a self-similar process with self-similarity parameter $H$. Its inverse hence transforms a self-similar process back to a stationary process. The inverse Lamperti transformation of a self-similar process $\{X(t)\}_{t\ge0}$ with self-similarity parameter $H\in(0,1)$ is thus given by
$$
Y(t)=e^{-Ht}X(e^{t}),~\mbox{for}~t\ge0.
$$
Since $\{Y(t)\}_{t\ge0}$ is stationary, the simulation of $\{X(t)\}_{t\ge0}$ can be obtained through the simulation of $\{Y(t)\}_{t\ge0}$, which is much simpler. This idea may provide hope for developing a general method for simulating self-similar processes. 

By using the Lamperti transformations, we suggest the following novel method to simulate a zero-mean Gaussian self-similar process $\{X(t)\}_{t\in[0,1]}$:
\begin{enumerate}
\item Given the length $n\ge1$, use Wood-Chan's method \cite{WoodChan1994} to simulate a sample path of the stationary Gaussian process
\begin{equation}
\label{def:U}
U\left(\frac{k}{n}\right)=n^{-H(k/n-1)}X(n^{k/n-1}),~\mbox{for}~k=1,\ldots,n
\end{equation}
via its autocovariance function
\begin{equation}
\label{def:gamma}
\gamma_U\left(\frac{k}{n}\right)=\mathbb Cov\left(U\left(\frac{1}{n}\right),U\left(\frac{k+1}{n}\right)\right).
\end{equation}
\item For $j=1,\ldots,n$, take
\begin{equation}
\label{def:Xtilde}
\widetilde X(\frac{j}{n}) =\left(\frac{j}{n}\right)^HU\left(\frac{\left\lfloor\left(\frac{\log(j/n)}{\log n}+1\right)n\right\rfloor}{n}\right),
\end{equation}
where $\lfloor\bullet\rfloor$ denotes the floor number. 
\end{enumerate}
We explain how the above algorithm works. (\ref{def:U}) is a modified version of the inverse Lamperti transformation of $\{X(t)\}_{t\in[0,1]}$, where the mapping $t\mapsto n^{t-1}$ is a bijection from $[0,1]$ to $[1/n,1]$. This mapping allows one to transform $(U(1/n),\ldots,U(1))$ back to $(X(1/n),\ldots,X(1))$. By (\ref{def:U}), we get
\begin{eqnarray}
\label{def:tildeX}
\left(\widetilde X(\frac{j}{n})\right)_{j=1,\ldots,n} &=&\left(\left(\frac{j}{n}\right)^HU\left(\frac{\left\lfloor\left(\frac{\log(j/n)}{\log n}+1\right)n\right\rfloor}{n}\right)\right)_{j=1,\ldots,n}\nonumber\\
&=& \left(\left(\frac{j}{n}\right)^Hn^{-H\left(\frac{\left\lfloor\left(\frac{\log(j/n)}{\log n}+1\right)n\right\rfloor}{n}-1\right)}X\left(n^{\frac{\left\lfloor\left(\frac{\log(j/n)}{\log n}+1\right)n\right\rfloor}{n}-1}\right)\right)_{j=1,\ldots,n}\nonumber\\
&=&\left(\left(\frac{j}{n}\right)^Hn^{-H\left(\frac{\left(\frac{\log(j/n)}{\log n}+1\right)n-\theta_{j,n}}{n}-1\right)}X\left(n^{\frac{\left(\frac{\log(j/n)}{\log n}+1\right)n-\theta_{j,n}}{n}-1}\right)\right)_{j=1,\ldots,n}\nonumber\\
&=&\left(n^{H\theta_{j,n}/n}X\left(\frac{j}{n}\cdot n^{-\theta_{j,n}/n}\right)\right)_{j=1,\ldots,n},
\end{eqnarray}
where 
\begin{equation}
\label{def:theta}
\theta_{j,n}=\left(\frac{\log(j/n)}{\log n}+1\right)n-\left\lfloor\left(\frac{\log(j/n)}{\log n}+1\right)n\right\rfloor\in(0,1)~ \mbox{for}~ j\in\{1,\ldots,n\}.
\end{equation}
Then $(\widetilde X(j/n))_{j=1,\ldots,n}$ are \enquote{good} approximations of $( X(j/n))_{j=1,\ldots,n}$ in the following senses: 
\begin{description}
\item[(i)] By the self-similarity, we see that 
\begin{equation}
\label{error1}
\widetilde X\left(\frac{j}{n}\right)\stackrel{\mbox{law}}{=} X\left(\frac{j}{n}\right)~ \mbox{for each}~ j\in\{1,\ldots,n\}. 
\end{equation}
This means that the marginal distributions of $(\widetilde X(j/n))_{j=1,\ldots,n}$ and $( X(j/n))_{j=1,\ldots,n}$ are equal. However, their joint distributions are in general distinct.\\
\item[(ii)] Assume that the sample paths of $\{X(t)\}_{t\in[0,1]}$ are almost surely H\"older continuous with H\"older's index $\beta\in(0,1)$, i.e., there is a positive random variable $C$ such that
\begin{equation}
\label{Holder}
\sup_{s,t\in[0,1]}|X(s)-X(t)|\le C|s-t|^{\beta}~\mbox{a.s.}.
\end{equation}
Then by the triangle inequality and (\ref{Holder}), for any $j\in\{1,\ldots,n\}$,
\begin{eqnarray*}
&&\left|n^{H\theta_{j,n}/n}X\left(\frac{j}{n}\cdot n^{-\theta_{j,n}/n}\right)-X\left(\frac{j}{n}\right)\right|\\
&&\le \left|n^{H\theta_{j,n}/n}-1\right|\left|X\left(\frac{j}{n}\cdot n^{-\theta_{j,n}/n}\right)\right|+\left|X\left(\frac{j}{n}\cdot n^{-\theta_{j,n}/n}\right)-X\left(\frac{j}{n}\right)\right|\\
&&\le \left(\sup_{s\in[0,1]}|X(s)|\right)\left|n^{H\theta_{j,n}/n}-1\right|+C\left(\frac{j}{n}\right)^{\beta-\varepsilon}\left| n^{-\theta_{j,n}/n}-1\right|^{\beta}.
\end{eqnarray*}
By the mean value theorem and the fact that $\lim_{n\to+\infty}n^{1/n}=1$, there is a constant $c_1>0$ which does not depend on $j$ nor on $n$, such that
$$
\left|n^{H\theta_{j,n}/n}-1\right|=\left|\left(n^{1/n}\right)^{H\theta_{j,n}}-\left(n^{1/n}\right)^0\right|= \frac{(n^{1/n})^{\eta_{j,n}}\log(n)}{n}H\theta_{j,n}\le c_1\frac{\log(n)}{n},
$$
where $\eta_{j,n}$ is some value in $(0,H\theta_{j,n})\subset(0,H)$ for any $n\ge1$ and $j\in\{1,\ldots,n\}$. In the same way, there is a constant $c_2>0$ which does not depend on $j$ nor on $n$, such that
$$
\left|n^{-\theta_{j,n}/n}-1\right|\le c_2\frac{\log(n)}{n}.
$$
It follows that, for any $j\in\{1,\ldots,n\}$,
\begin{equation}
\label{diff:X}
\left|n^{H\theta_{j,n}/n}X\left(\frac{j}{n}\cdot n^{-\theta_{j,n}/n}\right)-X\left(\frac{j}{n}\right)\right|\le \widetilde Cn^{-\beta}(\log(n))^{\beta},
\end{equation}
where the random variable 
$
\widetilde C=  c_1\sup_{s\in[0,1]}|X(s)|+c_2C
$ does not depend on $j$ nor on $n$.   Finally, (\ref{diff:X}) together with (\ref{def:tildeX}) yields
\begin{equation}
\label{X:approx}
\max_{j\in\{1,\ldots,n\}}\left|\widetilde X\left(\frac{j}{n}\right)-X\left(\frac{j}{n}\right)\right|\le \widetilde Cn^{-\beta}(\log(n))^{\beta}\xrightarrow[n\to+\infty]{a.s.}0.
\end{equation}
We know that $\{X(t)\}_{t\in[0,1]}$ is a Gaussian process. If all its sample paths are almost surely continuous, then by applying Dudley's theorem
 and Borell's inequality (see Pages 1445-1446 in \cite{Rosenbaum2008}, see also \cite{Ledoux2010}), we can show that any order moments of $\sup_{s\in[0,1]}|X(s)|$ and $C$ are finite. These facts together with (\ref{X:approx}) imply that
 \begin{equation}
 \label{error21}
 \mathbb E\left(\max_{j\in\{1,\ldots,n\}}\left|\widetilde X\left(\frac{j}{n}\right)-X\left(\frac{j}{n}\right)\right|\right)^2\le cn^{-2\beta}(\log(n))^{2\beta},
 \end{equation}
 where $c=\mathbb E|\widetilde C|^2<+\infty$.
\end{description}

\subsection{Application I: Simulation of Fractional Brownian Motion}
\label{sim:fBm}
The modified inverse Lamperti transformation approach can be applied to simulate sample paths of a fractional Brownian motion $\{B^H(t)\}_{t\in[0,1]}$, defined via (\ref{cov_fBm}). It suffices to note that
\begin{eqnarray}
\label{def:gamma_fBm}
\gamma_U\left(\frac{k}{n}\right)&=&\mathbb Cov\left(n^{-H(1/n-1)}B^H(n^{1/n-1}),n^{-H((k+1)/n-1)}B^H(n^{(k+1)/n-1})\right)\nonumber\\
&=&\frac{1}{2}n^{-H((k+2)/n-2)}\left(n^{2H(1/n-1)}+n^{2H((k+1)/n-1)}-\left|n^{1/n-1}-n^{(k+1)/n-1}\right|^{2H}\right)\nonumber\\
&=&\frac{1}{2}\left(n^{-Hk/n}+n^{-Hk/n}-|n^{-k/(2n)}-n^{-k/(2n)}|^{2H}\right),~\mbox{for}~k=1,\ldots,n
\end{eqnarray}
in the simulation algorithm (\ref{def:U}) - (\ref{def:Xtilde}). As the fractional Brownian motion is a Gaussian process having almost surely H\"older continuous sample paths with the H\"older index $H-\varepsilon$ for any $\varepsilon>0$ arbitrarily small, the error analysis results (\ref{error1}) - (\ref{error21}) hold with $\beta=H-\varepsilon$. We summarize the pseudocode for simulating $\{B^H(t)\}_{t\in[0,1]}$ below:\\
\begin{minipage}{16cm}
\begin{algorithm}[H] 
\label{algo:simulation_self_similar}
\caption{Simulation of a discrete sample path of the fractional Brownian motion  $(B^H(1/n),B^H(2/n),\ldots,B^H(1))$.}
\LinesNumbered 
\KwIn{\textbf{\textbf{\textit{the path length}} $n$; \textit{the Hurst parameter} $H$; \textit{the covariance function}} $\gamma_U$ given in (\ref{def:gamma_fBm}).}
\textbf{\textit{Use Wood-Chan's method \cite{WoodChan1994} to simulate a sample path of $U$ based on its autocovariance function $\gamma_U$:}\\
$U \longleftarrow \Big\{U(1/n),U(2/n)\ldots,U(1)\Big \}$}\;
\textbf{\textit{Simulate $(B^H(1/n),B^H(2/n),\ldots,B^H(1))$}:}\\
\For{$j = 1,\ldots,n$}{
$B^H(\frac{j}{n}) \longleftarrow\left(\frac{j}{n}\right)^HU\left(\frac{\left\lfloor\left(\frac{\log(j/n)}{\log n}+1\right)n\right\rfloor}{n}\right)$\;}
\KwOut{ A sample path $(B^H(1/n),B^H(2/n),\ldots,B^H(1))$.}
\end{algorithm}
\end{minipage}
Algorithm \ref{algo:simulation_self_similar} has been implemented in the Python library \enquote{fractal analysis} (\url{https://pypi.org/project/fractal-analysis/}) by Yujia Ding. The above Python library also implements \cite{balcerek2020testingfBm}'s approach for the hypothesis test of fractional Brownian motion.  The test tells whether the observed path is sampled from the fractional Brownian motion $\{B^H(t)\}_{t\in[0,1]}$ with some given $H$. Using this test, we compare our algorithm to Wood-Chan's method. We summarize the comparison results below: \\
\begin{table}[H]
\begin{center}
\begin{tabular}{ |l|c|c|c|c|  }
 \hline
 \multicolumn{5}{|c|}{Pass Rate of Balcerek and Burnecki's Hypothesis Test} \\
 \hline
 \hline
 \multirow{2}{*}{Hurst Parameter $H$} &
      \multicolumn{2}{|c|}{$\alpha=5\%$} &
      \multicolumn{2}{|c|}{$\alpha=1\%$}  \\
      \cline{2-5}
 & Wood-Chan &Algorithm \ref{algo:simulation_self_similar}& Wood-Chan &Algorithm \ref{algo:simulation_self_similar}\\
 \hline
 0.01   & 72.8\%    &97.8\%&   93.5\%&99.9\%\\
 0.2&   87.6\%  & 100\%  &97.0\%&100\%\\
 0.5 &90.3\% & 99.8\%&  97.5\%&100\%\\
 0.8    &91.8\% & 99.5\%&  98.2\%&100\%\\
 0.99&   93.8\%  & 99.4\%&98.6\%&100\%\\
 \hline
\end{tabular}
\caption{$\alpha$ is the significance level of the test. Each pass rate is obtained by testing $1000$ simulated sample paths of $(B^H(1/n),\ldots,B^H(1))$ with length $n=1024$.}
\label{table:comparison}
\end{center}
\end{table}
Although the hypothesis test has its own testing errors, Table \ref{table:comparison} shows that Algorithm \ref{algo:simulation_self_similar} outperforms Wood-Chan's method in view of the pass rate. We also observe that the performance of Wood-Chan's method is not consistent with respect to the parameter $H$: greater is the value of $H$, greater is the pass rate. This is because Wood-Chan's method is based on simulating the increments of fractional Brownian motion, $\{B^H((j+1)/n)-B^H(j/n)\}_{j=0,\ldots,n-1}$, which have higher variance when $H$ is smaller. However, Algorithm \ref{algo:simulation_self_similar} does not have such issue. Below we illustrate some sample paths simulated by Algorithm \ref{algo:simulation_self_similar}.\\

\begin{figure}[ht!]  
   ~~~~~~~~~~~~~~~~~~~~~~$B^H(t)$ with $H=0.2$~~~~~~~~~~~~~~~~~~~~~~~~~~~~~~~~~~~~~~~~~~~~~~~~~~~~~~$B^H(t)$ with $H=0.8$\\
    \includegraphics[width=.45\textwidth]{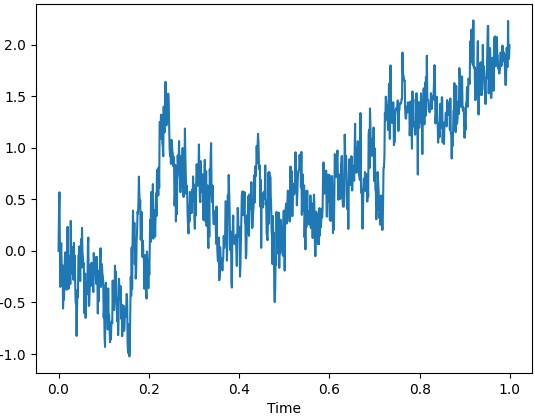}
         \hfill
   \includegraphics[width=.45\textwidth]{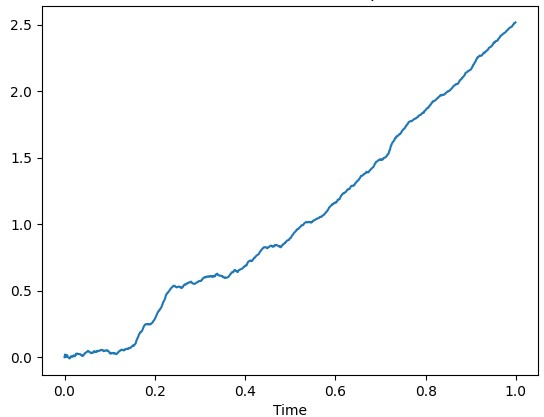}
   
   \caption{Sample paths of $\{B^H(t)\}_{t\in[0,1]}$ simulated by using Algorithm \ref{algo:simulation_self_similar}.}
   \label{fig:fBm}
        \end{figure}
        
From Figure \ref{fig:fBm}, we see that the H\"older's regularity of the simulated sample paths indeed changes via the value of $H$. 

\subsection{Application II: Simulation of Sub-fractional Brownian Motion}
\label{sim:subfBm}
Unlike fractional Brownian motion, the sub-fractional Brownian motion $\{S^H(t)\}_{t\in[0,1]}$ defined in (\ref{def:sfBm}) has non-stationary increments, which makes Wood-Chan's method no longer work. However, the modified inverse Lamperti transformation method applies here. In the algorithm (\ref{def:U}) - (\ref{def:Xtilde}), we take
\begin{eqnarray}
\label{def:gamma_sfBm}
\gamma_U\left(\frac{k}{n}\right)&=&\mathbb Cov\left(n^{-H(1/n-1)}S^H(n^{1/n-1}),n^{-H((k+1)/n-1)}S^H(n^{(k+1)/n-1})\right)\nonumber\\
&=&n^{-H((k+2)/n-2)}\left(n^{2H(1/n-1)}+n^{2H((k+1)/n-1)}\right.\nonumber\\
&&\quad\left.-\frac{1}{2}\left(\left(n^{1/n-1}+n^{(k+1)/n-1}\right)^{2H}+\left|n^{1/n-1}-n^{(k+1)/n-1}\right|^{2H}\right)\right)\nonumber\\
&=&n^{-Hk/n}+n^{Hk/n}-\frac{1}{2}\left((n^{-k/(2n)}+n^{k/(2n)})^{2H}+|n^{-k/(2n)}-n^{k/(2n)}|^{2H}\right),\nonumber\\
\end{eqnarray}
for $k=1,\ldots,n$. 
From the sub-fractional Brownian motion's covariance structure and the Kolmogorov-\v{C}entsov theorem, we know that $\{S^H(t)\}_{t\in[0,1]}$ is also a Gaussian process whose sample paths are almost surely $(H-\varepsilon)$-H\"older continuous for arbitrarily small $\varepsilon>0$. Therefore, the error analysis results (\ref{error1}) - (\ref{error21}) hold with $\beta=H-\varepsilon$ for the sub-fractional Brownian motion. We summarize the pseudocode for simulating $\{S^H(t)\}_{t\in[0,1]}$ below:\\
\begin{minipage}{16cm}
\begin{algorithm}[H] 
\label{algo:simulation_sub_self_similar}
\caption{Simulation of a discrete sample path of the sub-fractional Brownian motion  $(S^H(1/n),S^H(2/n),\ldots,S^H(1))$.}
\LinesNumbered 
\KwIn{\textbf{\textbf{\textit{the path length}} $n$; \textit{the self-similarity index} $H$; \textit{the covariance function}} $\gamma_U$ given in (\ref{def:gamma_sfBm}).}
\textbf{\textit{Use Wood-Chan's method  to simulate a sample path of $U$ based on its autocovariance function $\gamma_U$:}\\
$U \longleftarrow \Big\{U(1/n),U(2/n)\ldots,U(1)\Big \}$}\;
\textbf{\textit{Simulate $(S^H(1/n),S^H(2/n),\ldots,S^H(1))$}:}\\
\For{$j = 1,\ldots,n$}{
$S^H(\frac{j}{n}) \longleftarrow\left(\frac{j}{n}\right)^HU\left(\frac{\left\lfloor\left(\frac{\log(j/n)}{\log n}+1\right)n\right\rfloor}{n}\right)$\;}
\KwOut{ A sample path $(S^H(1/n),S^H(2/n),\ldots,S^H(1))$.}
\end{algorithm}
\end{minipage}
Algorithm \ref{algo:simulation_sub_self_similar} has also been implemented in the Python library \enquote{fractal analysis}. Since the test of sub-fractional Brownian motion has not yet been found in the literature, we only illustrate some sample paths simulated by Algorithm \ref{algo:simulation_sub_self_similar} below.\\

\begin{figure}[H]  
   ~~~~~~~~~~~~~~~~~~~~~~$S^H(t)$ with $H=0.2$~~~~~~~~~~~~~~~~~~~~~~~~~~~~~~~~~~~~~~~~~~~~~~~~~~~~~~$S^H(t)$ with $H=0.8$\\
    \includegraphics[width=.45\textwidth]{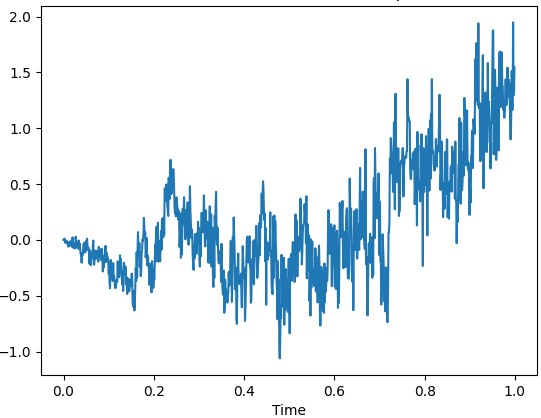}
         \hfill
   \includegraphics[width=.45\textwidth]{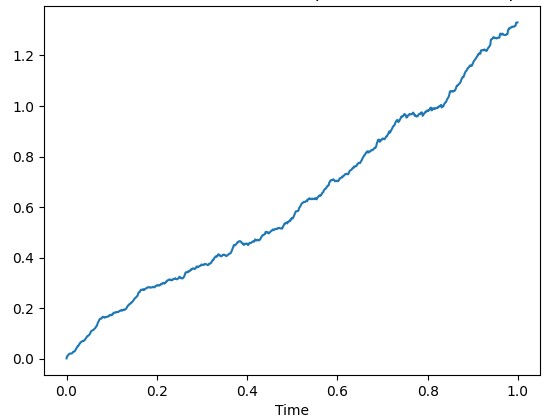}
   
   \caption{Sample paths of $\{S^H(t)\}_{t\in[0,1]}$ simulated by using Algorithm \ref{algo:simulation_sub_self_similar}.}
   \label{fig:sfBm}
        \end{figure}
        
Figure \ref{fig:sfBm} supports the fact that the sub-fractional Brownian motion is self-similar with the self-similarity index $H$.
        
 \section{Conclusion}
 \label{Sec7}
 There is a growing body of research on various self-similar processes, including the introduction of new types of self-similar processes, studying their paths features, and developing statistical inference tools. There is not yet a promising methodology that can efficiently simulate their sample paths. For example, although bi-fractional Brownian motion and tri-fractional Brownian motion have been receiving growing attention, simulation methods have yet to be developed for them. The first goal of this paper then has been to illuminate a direction of research that aims to shed light on the development of general simulation strategies for the above self-similar processes. The prototypical examples we have chosen to investigate are fractional Brownian motion, sub-fractional Brownian motion, multifractional Brownian motion, and linear fractional stable motion. Based on the above investigation of the existing simulation approaches, we believe that working on the transformation of self-similar processes to stationary processes may point in the direction of an efficient generating algorithm of self-similar processes. Therefore, the Lamperti transformation is a potential powerful tool to be used to simulate a general self-similar process. We then have proposed a new method to simulate self-similar processes. This method uses a modified version of the inverse Lamperti transformation to transform the self-similar process to a stationary Gaussian process. Simulating the self-similar process is then equivalent to simulate the latter stationary process. The second contribution of this paper has been to successfully apply this novel method to simulate the fractional Brownian motion and sub-fractional Brownian motion. Comparison results show that our new method outperforms Wood-Chan's method via the hypothesis test pass rate. We believe that this idea can be further extended to simulate a large class of self-similar processes which are Gaussian processes with non-stationary increments or are non-Gaussian processes, including bi-fractional Brownian motion, tri-fractional Brownian motion, and linear fractional stable motion.

\section*{\textit{Acknowledgment:}}

\textit{The authors thank Yujia Ding for helping to implement the simulation and testing algorithms in the Python library.}

\bibliographystyle{plain}
\bibliography{ref}

@article{elsaesser2010observed,
  title={Observed self-similarity of precipitation regimes over the tropical oceans},
  author={Elsaesser, Gregory S and Kummerow, Christian D and L’Ecuyer, Tristan S and Takayabu, Yukari N and Shige, Shoichi},
  journal={Journal of Climate},
  volume={23},
  number={10},
  pages={2686--2698},
  year={2010}
}

@article{dobric2006fractional,
  title={Fractional {B}rownian fields, duality, and martingales},
  author={Dobri{\'c}, Vladimir and Ojeda, Francisco M},
  journal={Lecture Notes-Monograph Series},
  pages={77--95},
  year={2006},
  publisher={JSTOR}
}

@article{dobrushin1979non,
  title={Non-central limit theorems for non-linear functional of {G}aussian fields},
  author={Dobrushin, Roland Lvovich and Major, P{\'e}ter},
  journal={Zeitschrift f{\"u}r Wahrscheinlichkeitstheorie und verwandte Gebiete},
  volume={50},
  pages={27--52},
  year={1979},
  publisher={Springer}
}

@article{taqqu1979convergence,
  title={Convergence of integrated processes of arbitrary Hermite rank},
  author={Taqqu, Murad S},
  journal={Zeitschrift f{\"u}r Wahrscheinlichkeitstheorie und verwandte Gebiete},
  volume={50},
  number={1},
  pages={53--83},
  year={1979},
  publisher={Springer}
}

@article{laskin2003fractional,
  title={Fractional {P}oisson process},
  author={Laskin, Nick},
  journal={Communications in Nonlinear Science and Numerical Simulation},
  volume={8},
  number={3-4},
  pages={201--213},
  year={2003},
  publisher={Elsevier}
}

@article{sun2018self,
  title={Self-similarity and the dynamics of coarsening in materials},
  author={Sun, Yue and Andrews, W Beck and Thornton, Katsuyo and Voorhees, Peter W},
  journal={Scientific Reports},
  volume={8},
  number={1},
  pages={17940},
  year={2018},
  publisher={Nature Publishing Group UK London}
}

@article{haltas2011scale,
  title={Scale invariance and self-similarity in hydrologic processes in space and time},
  author={Haltas, I and Kavvas, ML},
  journal={Journal of Hydrologic Engineering},
  volume={16},
  number={1},
  pages={51--63},
  year={2011},
  publisher={American Society of Civil Engineers}
}

@article{rostek2013note,
  title={A note on the use of fractional {B}rownian motion for financial modeling},
  author={Rostek, Stefan and Sch{\"o}bel, R},
  journal={Economic Modelling},
  volume={30},
  pages={30--35},
  year={2013},
  publisher={Elsevier}
}

@article{ma2013schoenberg,
  title={The {S}choenberg--{L}e\'evy kernel and relationships among fractional {B}rownian motion, bifractional {B}rownian motion, and others},
  author={Ma, Chungsheng},
  journal={Theory of Probability \& Its Applications},
  volume={57},
  number={4},
  pages={619--632},
  year={2013},
  publisher={SIAM}
}

@phdthesis{jeong2002modelling,
 title={Modelling of Self-similar Teletraffic for Simulation},
  author={Jeong, Hae-Duck Joshua},
  year={2002},
  school={University of Canterbury}
}

@article{stoev2004simulation,
  title={Simulation methods for linear fractional stable motion and {FARIMA} using the Fast Fourier Transform},
  author={Stoev, Stilian and Taqqu, Murad S},
  journal={Fractals},
  volume={12},
  number={01},
  pages={95--121},
  year={2004},
  publisher={World Scientific}
}

@article{bierme2008fourier,
  title={Fourier series approximation of linear fractional stable motion},
  author={Bierm{\'e}, Hermine and Scheffler, Hans-Peter},
  journal={Journal of Fourier Analysis and Applications},
  volume={14},
  number={2},
  pages={180--202},
  year={2008},
  publisher={Springer}
}

@article{wu2004simulating,
  title={Simulating sample paths of linear fractional stable motion},
  author={Wu, Wei Biao and Michailidis, George and Zhang, Danlu},
  journal={IEEE Transactions on Information Theory},
  volume={50},
  number={6},
  pages={1086--1096},
  year={2004},
  publisher={IEEE}
}

@inproceedings{morozewicz2015simulation,
  title={On the simulation of sub-fractional {B}rownian motion},
  author={Morozewicz, Aneta and Filatova, Darya},
  booktitle={2015 20th International Conference on Methods and Models in Automation and Robotics (MMAR)},
  pages={400--405},
  year={2015},
  organization={IEEE}
}

@article{kuang2017asymptotic,
  title={Asymptotic behavior of weighted cubic variation of sub-fractional {B}rownian motion},
  author={Kuang, Nenghui and Xie, Huantian},
  journal={Communications in Statistics-Simulation and Computation},
  volume={46},
  number={1},
  pages={215--229},
  year={2017},
  publisher={Taylor \& Francis}
}

@inproceedings{ayache2011multiparameter,
  title={Multiparameter multifractional {B}rownian motion: local nondeterminism and joint continuity of the local times},
  author={Ayache, Antoine and Shieh, Narn-Rueih and Xiao, Yimin},
  booktitle={Annales de l'IHP Probabilit{\'e}s et Statistiques},
  volume={47},
  number={4},
  pages={1029--1054},
  year={2011}
}

@book{embrechts2009selfsimilar,
  title={Selfsimilar Processes},
  author={Embrechts, Paul},
  year={2009},
  publisher={Princeton University Press}
}

@book{beran2017statistics,
  title={Statistics for Long-memory Processes},
  author={Beran, Jan},
  year={2017},
  publisher={Routledge}
}

@article{davies1987tests,
  title={Tests for {H}urst effect},
  author={Davies, Robert B and Harte, D.S.},
  journal={Biometrika},
  volume={74},
  number={1},
  pages={95--101},
  year={1987},
  publisher={Oxford University Press}
}

@phdthesis{michna2000ruin,
  title={Ruin Probabilities and First Passage Times for Self-similar Processes},
  author={Michna, Zbigniew},
  year={2000},
  school={Department of Mathematical Statistics, Lund University}
}

@article{michna1998self,
  title={Self-similar processes in collective risk theory},
  author={Michna, Zbigniew},
  journal={International Journal of Stochastic Analysis},
  volume={11},
  number={4},
  pages={429--448},
  year={1998},
  publisher={Wiley Online Library}
}

@article{michna1999tail,
  title={On tail probabilities and first passage times for fractional {B}rownian motion},
  author={Michna, Zbigniew},
  journal={Mathematical Methods of Operations Research},
  volume={49},
  pages={335--354},
  year={1999},
  publisher={Springer}
}

@article{ding2023linear,
  title={Linear multifractional stable sheets in the broad sense: existence and joint continuity of local times},
  author={Ding, Yujia and Peng, Qidi and Xiao, Yimin},
  journal={Bernoulli},
  volume={29},
  number={1},
  pages={785--814},
  year={2023},
  publisher={Bernoulli Society for Mathematical Statistics and Probability}
}

@article{ayache2007local,
  title={Local and asymptotic properties of linear fractional stable sheets},
  author={Ayache, Antoine and Roueff, Fran{\c{c}}ois and Xiao, Yimin},
  journal={Comptes rendus. Math{\'e}matique},
  volume={344},
  number={6},
  pages={389--394},
  year={2007}
}

@article{ayache2007wavelet,
  title={Wavelet construction of generalized multifractional processes},
  author={Ayache, Antoine and Jaffard, St{\'e}phane and Taqqu, Murad S},
  journal={Revista Matematica Iberoamericana},
  volume={23},
  number={1},
  pages={327--370},
  year={2007}
}

@book{asmussen1998stochastic,
  title={Stochastic Simulation with A View towards Stochastic Processes},
  author={Asmussen, S{\o}ren},
  year={1998},
  publisher={University of Aarhus. Centre for Mathematical Physics and Stochastics}
}

@inproceedings{chan1998simulation,
  title={Simulation of multifractional {B}rownian motion},
  author={Chan, Grace and Wood, Andrew T.A.},
  booktitle={COMPSTAT},
  pages={233--238},
  year={1998},
  organization={Springer}
}

@article{lamperti1962semi,
  title={Semi-stable stochastic processes},
  author={Lamperti, John},
  journal={Transactions of the American Mathematical Society},
  volume={104},
  number={1},
  pages={62--78},
  year={1962},
  publisher={JSTOR}
}

@article {Benassi1997Elliptic,
    AUTHOR = {Benassi, Albert and Jaffard, St\'{e}phane and Roux, Daniel},
     TITLE = {Elliptic {G}aussian random processes},
  JOURNAL = {Revista Matem\'{a}tica Iberoamericana},
    VOLUME = {13},
      YEAR = {1997},
    NUMBER = {1},
     PAGES = {19--90}
}

@article{bojdecki2004sub,
  title={Sub-fractional {B}rownian motion and its relation to occupation times},
  author={Bojdecki, Tomasz and Gorostiza, Luis G and Talarczyk, Anna},
  journal={Statistics \& Probability Letters},
  volume={69},
  number={4},
  pages={405--419},
  year={2004},
  publisher={Elsevier}
}

@article{houdre2003example,
  title={An example of infinite dimensional quasi-helix},
  author={Houdr{\'e}, Christian and Villa, Jos{\'e}},
  journal={Contemporary Mathematics},
  volume={336},
  pages={195--202},
  year={2003},
  publisher={Citeseer}
}

@article{vehel1998introduction,
  title={Introduction to the multifractal analysis of images},
  author={ L{\'e}vy-V{\'e}hel, Jacques},
  journal={Fractal Image Encoding and Analysis},
  volume={159},
  pages={299--341},
  year={1998},
  publisher={Springer}
}

@article{kolmogorov1940wienersche,
  title={Wienersche spiralen und einige andere interessante kurven in hilbertscen raum, cr (doklady)},
  author={Kolmogorov, Andrei N},
  journal={Acad. Sci. URSS (NS)},
  volume={26},
  pages={115--118},
  year={1940}
}

@article{Mandelbrot1968,
  title={Fractional {B}rownian motions, fractional noises and applications},
  author={Mandelbrot, B and van Ness, J W},
  journal={SIAM Review},
  volume={10},
  number={4},
  pages={422--437},
  year={1968}
}

@article{Rosenbaum2008,
  title={Estimation of the volatility persistence in a discretely observed diffusion model},
  author={Rosenbaum, Mathieu},
  journal={Stochastic Processes and their Applications},
  volume={118},
  number={8},
  pages={1434--1462},
  year={2008},
  publisher={Elsevier}
}

@article{coeurjolly2013wavelet,
  title={Wavelet analysis of the multivariate fractional {B}rownian motion},
  author={Coeurjolly, Jean-Fran{\c{c}}ois and Amblard, Pierre-Olivier and Achard, Sophie},
  journal={ESAIM: Probability and Statistics},
  volume={17},
  pages={592--604},
  year={2013},
  publisher={EDP Sciences}
}

@article{jeong2003generation,
  title={Generation of self-similar processes for simulation studies of telecommunication networks},
  author={Jeong, Hae-Duck J and Pawlikowski, K and McNickle, DC},
  journal={Mathematical and Computer Modelling},
  volume={38},
  number={11-13},
  pages={1249--1257},
  year={2003},
  publisher={Elsevier}
}

@inproceedings{jeong1999fast,
  title={Fast self-similar teletraffic generation based on {FGN} and wavelets},
  author={Jeong, H-DJ and McNickle, Don and Pawlikowski, Krzysztof},
  booktitle={IEEE International Conference on Networks. ICON'99 Proceedings (Cat. No. PR00243)},
  pages={75--82},
  year={1999},
  organization={IEEE}
}

@book{doukhan2002theory,
  title={Theory and Applications of Long-range Dependence},
  author={Doukhan, Paul and Oppenheim, George and Taqqu, Murad},
  year={2002},
  publisher={Springer Science \& Business Media}
}

@article{albeverio2012fractional,
  title={On fractional {B}rownian motion and wavelets},
  author={Albeverio, S and Jorgensen, Palle ET and Paolucci, Anna Maria},
  journal={Complex Analysis and Operator Theory},
  volume={6},
  number={1},
  pages={33--63},
  year={2012},
  publisher={Springer}
}

@article{flandrin1992wavelet,
  title={Wavelet analysis and synthesis of fractional {B}rownian motion},
  author={Flandrin, Patrick},
  journal={IEEE Transactions on Information Theory},
  volume={38},
  number={2},
  pages={910--917},
  year={1992},
  publisher={IEEE}
}

@article{whitcher2001simulating,
  title={Simulating {G}aussian stationary processes with unbounded spectra},
  author={Whitcher, Brandon},
  journal={Journal of Computational and Graphical Statistics},
  volume={10},
  number={1},
  pages={112--134},
  year={2001},
  publisher={Taylor \& Francis}
}

@article{mallat1989theory,
  title={A theory for multiresolution signal decomposition: the wavelet representation},
  author={Mallat, Stephane G},
  journal={IEEE Transactions on Pattern Analysis and Machine Intelligence},
  volume={11},
  number={7},
  pages={674--693},
  year={1989},
  publisher={Ieee}
}

@article{norros1999elementary,
  title={An elementary approach to a {G}irsanov formula and other analytical results on fractional {B}rownian motions},
  author={Norros, Ilkka and Valkeila, Esko and Virtamo, Jorma},
journal={Bernoulli},
volume={5},
  number={4},
  pages={571--587},
  year={1999}
}

@book{taqqu1979self,
  title={Self-similar processes and related ultraviolet and infrared catastrophes},
  author={Taqqu, Murad S},
  year={1979},
  publisher={Technical Report - School of Operations Research and Industrial Engineering, College of Engineering, Cornell University}
}

@article{WoodChan1994,
  title={Simulation of stationary {G}aussian processes in $[0, 1]^d$},
  author={Wood, Andrew T A and Chan, Grace},
  journal={Journal of Computational and Graphical Statistics},
  volume={3},
  number={4},
  pages={409--432},
  year={1994},
  publisher={Taylor \& Francis Group}
}

@article{Stoev2006,
  title={How rich is the class of multifractional {B}rownian motions?},
  author={Stoev, Stilian A and Taqqu, Murad S},
  journal={Stochastic Processes and their Applications},
  volume={116},
  number={2},
  pages={200--221},
  year={2006},
  publisher={Elsevier}
}

@article{Jin2018,
author = {Sixian Jin and Qidi Peng and Henry Schellhorn},
year = {2018},
pages = {113-140},
title = {Estimation of the pointwise {H}\"older exponent of hidden multifractional {B}rownian motion using wavelet coefficients},
volume = {21},
number = {1},
journal = {Statistical Inference for Stochastic Processes}
}

@article{balcerek2020testingfBm,
  title={Testing of fractional {B}rownian motion in a noisy environment},
  author={Balcerek, Micha{\l} and Burnecki, Krzysztof},
  journal={Chaos, Solitons \& Fractals},
  volume={140},
  pages={110097},
  year={2020},
  publisher={Elsevier}
}

@book{Taqqu1994,
  author = {G. Samorodnitsky and M. S. Taqqu},
  year = {1994},
  title = {Stable {N}on-{G}aussian {R}andom {P}rocesses},
  publisher = {Chapman \& Hall/CRC}
}

@book{Ledoux2010,
  title={Probability in Banach Spaces: Isoperimetry and Processes},
  author={Ledoux, Michel and Talagrand, Michel},
  year={2013},
  publisher={Springer Science \& Business Media}
}

@article{Boufoussi2008,
  title		= {Path properties of a class of locally asymptotically
self-similar processes},
  author	= {Boufoussi, Brahim  and Dozzi, Marco and  Guerbaz, Raby},
  journal	= {Electronic Journal of Probability},
  volume	= {13},
  number	= {29},
  pages		= {898--921},
  year		= {2008}
}

\end{document}